\documentstyle[12pt,amssymb,amscd,psfig]{amsart}
\newcommand{\summ}{\sum\limits}

\newtheorem{thm}{Theorem}
\newtheorem{corol}[thm]{Corollary}
\newtheorem{lemma}[thm]{Lemma}

\theoremstyle{definition}

\theoremstyle{remark}

\newcommand{\co}{\colon\thinspace}    %  Colon with correct spacing for maps.
\newcommand{\sphere}{
$$\begin{picture}(240,120)  \footnotesize
    \put(-35,5)       {\psfig{figure=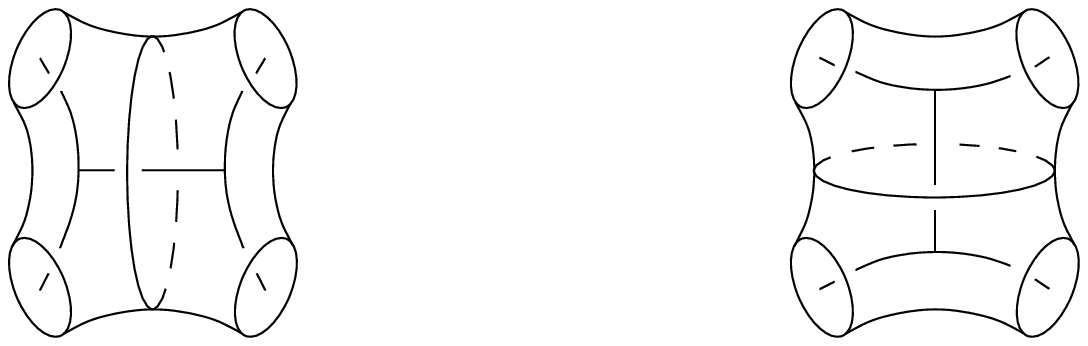}}
     \put(-16,72)      {$a_2$}
     \put(19,72)       {$a_3$}
     \put(-16,32)      {$a_1$}
     \put(20,32)       {$a_4$}
     \put(5,56)        {$b$}
     \put(58,50)       {{\normalsize $ =\summ_{b'}^{ }\left\{ \begin{array}{ccc} a_1 & a_2 & b'\\ a_3 & a_4 & b \end{array} \right\} $ }}
     \put(210,73)       {$a_2$}
     \put(248,72)       {$a_3$}
     \put(209,31)       {$a_1$}
     \put(248,31)       {$a_4$}
     \put(235,62)       {$b'$}
\end{picture}$$}

\title[The asymptotics of quantum $SU(2)$ representations]{On the
asymptotics of quantum ${\mathbf SU(2)}$ representations of mapping
class groups.}
\author[Michael Freedman and Vyacheslav Krushkal]{Michael Freedman$^{*}$ and Vyacheslav Krushkal$^{\dagger}$}
\thanks{$^{*}$ Microsoft Research, Redmond WA 98052}
\thanks{$^{\dagger}$ University of Virginia, Charlottesville VA 22904-4137}
\thanks{$^{\dagger}$ Partially supported by NSF grant DMS-0306934}
\begin{document}

\maketitle

\noindent
{\bf Abstract.} We investigate the rigidity and asymptotic
properties of quantum $SU(2)$ representations of mapping class
groups. In the spherical braid group case the trivial
representation is not isolated in the family of quantum $SU(2)$
representations. In particular, they may be used to give an explicit check
that spherical braid groups and hyperelliptic mapping class groups do not
have Kazhdan's property (T).
On the other hand,
the representations of the mapping class group of the
torus do not have almost invariant vectors, in fact
they converge to the metaplectic representation of $SL(2, {\mathbb
Z})$ on $L^2({\mathbb R})$. As a consequence we obtain a curious
analytic fact about the Fourier transform on ${\mathbb R}$ which may
not have been previously observed.

\medskip

\section{Introduction and statement of results}
The $(2+1)$-dimensional topological quantum field theories provide a
family of interesting projective representations of mapping class
groups of surfaces. Some properties of these representations have
been established: it is known that the $SU(2)$ quantum
representations at prime roots of unity are irreducible \cite{R},
and are asymptotically faithful \cite{A}, \cite{FWW}. However
overall these representations remain relatively unexplored. In this
paper we study their asymptotic properties from a different
perspective, focusing in particular on rigidity of representations
of mapping class groups.

There is a well established analogy between mapping class groups and
arithmetic groups, cf \cite{I}. This analogy involves both the
algebraic properties such as the Tits alternative and the finiteness
of the virtual cohomological dimension and the geometric aspects:
the action of the mapping class groups on Teichm\"{u}ller spaces
compared to the action of lattices in semisimple Lie groups on the
symmetric spaces. While some of the properties of mapping class
groups are similar to those of higher rank lattices, some are
analogous to lattices of rank one, cf \cite{FLM}.

In this context it is an interesting open question whether mapping
class groups have Kazhdan's property (T). In the genus one case the
group is $SL(2,{\mathbb Z})$ and doesn't have property (T). The
affirmative answer for genus greater than one would be similar to
the case of higher rank lattices (such as $SL(n,{\mathbb Z})$,
$n\geq 2$) which have property (T) \cite{K}. We propose an approach
to this question based on quantum representations of mapping class
groups. (See section 6.)

{\bf Notation}.
Given a surface $S$ of genus $g$ with $b$ boundary components and $n$
marked points, its mapping class group ${\cal M}^b_{g,n}$ is the group
of orientation-preserving diffeomorhisms of $S$ fixing $\partial S$ pointwise
and preserving the set of marked points, modulo isotopy fixing
$\partial S$ and marked points.

\smallskip

\begin{thm} \label{property T theorem} \sl
The quantum $SU(2)$ representations of ${\cal M}_0^n$ have almost 
invariant vectors. In particular, they may be used to explicitly check
that the spherical mapping class groups ${\cal M}_{0,n}$ and the 
hyperelliptic mapping class groups do not have Kazhdan's property (T).
\end{thm}

\smallskip

In the case $n=4$ this implies that ${\cal M}_2$ doesn't have property (T).
The fact that the spherical braid groups and the hyperelliptic mapping 
class groups do not have property (T) may also be 
proved directly, cf \cite{Ko}, using the fact that $SL(2,{\mathbb Z})$ does not 
have it. A more detailed
discussion of property (T) and the proof of theorem \ref{property T theorem}
are given in section \ref{property T}.

$PSL(2,{\mathbb Z})$ occurs as a mapping class group in two guises:
as ${\cal M}_{0,4}$ and as ${\cal M}_{1,0}$ modulo center. Therefore the quantum
representations give rise to two families of representations of
$PSL(2,{\mathbb Z})$. Analyzing the representations of ${\cal M}_{0,4}$
(and more generally of ${\cal M}_{0,n}$) yields the proof of theorem \ref{property T theorem} above.
In contrast, the representations of ${\cal M}_{1,0}$ exhibit
a rather different behavior.

\smallskip

\begin{thm} \label{convergence theorem} \sl
The quantum $SU(2)$ representations of the mapping class group of
the closed torus (isomorphic to $SL(2,{\mathbb Z})$) converge, as
the level $k\longrightarrow\infty$, to (an irreducible component of)
the metaplectic representation of $SL(2,{\mathbb Z})$ on
$L^2({\mathbb R})$.
\end{thm}

\smallskip

Both the quantum representations and the metaplectic representation are projective, however they may be modified to be genuine linear representations
of $SL(2,{\mathbb Z})$ (and actually of $PSL(2,{\mathbb Z})$.)
The convergence may be understood in both contexts: for representations of
the central extension of $SL(2,{\mathbb Z})$ by ${\mathbb Z}$, and for modified
representations of $SL(2,{\mathbb Z})$.
Convergence is in the Fell topology on the space of representations (see below).

Comparing the quantum $SU(2)$ representations of ${\cal M}_{1,0}$ with
the metaplectic representation over finite fields yields the following result.

\smallskip

\begin{lemma} \label{factor} \sl
For prime $p$ the quantum $SU(2)$ representation ${\rho}_p$ of ${\cal M}_{1,0}
\cong SL(2,{\mathbb Z})$ factors through $SL(2,{\mathbb Z}/p)$.
\end{lemma}

\smallskip

The proof shows that ${\rho}_p$ equals the metaplectic representation
over ${\mathbb Z}/p{\mathbb Z}$ composed with the mod $p$ reduction
$SL(2,{\mathbb Z})\longrightarrow SL(2,{\mathbb Z}/p{\mathbb Z})$.
This fact has been independently observed by Larsen and Wang \cite{LW}.
In particular this gives a new proof of the result \cite{G}
that the image of the representations ${\rho}_p$ is finite.

A corollary, to be contrasted with theorem \ref{property T theorem},
is that the quantum $SU(2)$ representations ${\rho}_k$ of ${\cal M}_{1,0}$
do not have almost invariant vectors.

\smallskip

\begin{corol} \label{no almost invariant vectors} \sl
The trivial representation of $SL(2,{\mathbb Z})$ is isolated in the family
of quantum $SU(2)$ representations of ${\cal M}_{1,0}$.
\end{corol}

\smallskip

Another application concerns Fourier analysis:
we show that there exists a subset of the real line of constant density
such that no set similar to it supports functions almost invariant under
the Fourier transform. The proofs of this result (stated as
theorem \ref{analytic} in section \ref{tau}) and of corollary
\ref{no almost invariant vectors} follow from property $\tau$
(a close relative of Kazhdan's property T) for $SL(2,{\mathbb Z})$.
This discussion and the proof of lemma \ref{factor} are contained in
section \ref{tau}.

The rest of the paper is organized as follows.
The aspects of quantum $SU(2)$ representations of mapping class groups,
relevant to our applications, are summarized in section \ref{quantum}.
Section \ref{property T} recalls the definition and basic properties of
Kazhdan's property (T) and Fell topology on the space of representations,
and gives the proof of theorem \ref{property T theorem}.
The background on metaplectic representation and the proof of theorem
\ref{convergence theorem} are given in section \ref{metaplectic}.
The paper concludes with a few comments and open questions.

\smallskip

\noindent
{\bf Acknowledgements}. We would like to thank Peter Sarnak for directing our
attention to the metaplectic representation.

The second author is grateful to the Theory Group
at Microsoft Research for warm hospitality and support during the work on this
project.

\smallskip

%\vfill

%\pagebreak

\section{Quantum $SU(2)$ representations} \label{quantum}

The purpose of this section is to introduce notations and to
describe basic properties of quantum $SU(2)$ representations
necessary for our applications.
We refer the reader to \cite{BHMV}, \cite{Tu}, \cite{Wa}
for a detailed introduction to TQFT's.

The $SU(2)$ level $k$ representation is denoted ${\rho}_k$.
Set $r=k+2$ and define $A=e^{2{\pi}i/4r}$.
Let ${\Sigma}_g$ be a closed surface of genus $g$ and let $H$
be a handlebody with $\partial H={\Sigma}_g$.
Fix a trivalent spine of $H$, then the basis vectors are given
by edge labellings of the spine, where each label $j$ satisfies
$0\leq j\leq r-2$, with the admissibility conditions on labels $a$,
$b$, $c$ at each vertex:

(a) the labels $a$, $b$ and $c$ satisfy the three triangle
inequalities,\\
(b) $a+b+c$ is even,\\
(c) the ``quantum cut off'' $a+b+c\leq 2r-4.$

If $\Sigma$ is not closed choose a labelling of the boundary components
by integers $0\leq l_{\partial}\leq r-2$. Cap off each boundary
component of $\Sigma$ by a disk, denoting the resulting closed surface
by $\widehat{\Sigma}$, and let $H$ be a handlebody bounded by
$\widehat{\Sigma}$. The basis vectors for $\Sigma$ are defined
analogously to the closed case above with the additional requirement
that the spine of $H$ meets each cap of $\widehat{\Sigma}$ once
and the corresponding labellings agree with the fixed boundary
labels $l_{\partial}$.

Let $\alpha$ be a curve in ${\Sigma}_g$ which bounds in $H$. Assume the disk
bounded by $\alpha$ intersects the spine in a point, and the label
of the corresponding edge is $l$. Then the basis vectors
are eigenvectors for the action of the Dehn twist ${\tau}_{{\alpha}}$
with the corresponding eigenvalue equal to
\begin{equation} \label{twist}
A^{-l(l+2)}=e^{-{\pi}i l(l+2)/2r}.
\end{equation}

In the genus one case
${\cal M}_1\cong SL(2,{\mathbb Z})$ is generated by the
Dehn twists ${\tau}_{\alpha}$, ${\tau}_{\beta}$ where $\alpha$, $\beta$
form a standard symplectic basis of curves on the torus.
Consider the solid torus such that $\alpha$ bounds and take
the spine to be the core of the solid torus, so there is just one label
$l$, $0\leq l\leq r-2$. The Dehn twist ${\tau}_{\alpha}$ acts according to the
formula (\ref{twist}) above. The
Dehn twist ${\tau}_{\beta}$ is the conjugate of ${\tau}_{\alpha}$ by the
$S$-matrix corresponding to the homeomorphism of the torus exchanging
${\alpha}$ and ${\beta}$:
$\left( \begin{smallmatrix} 0 & 1\\ -1 & 0 \end{smallmatrix} \right).$
The $S$-matrix equals the imaginary part of the discrete Fourier transform
\cite{Wi}:
\begin{equation} \label{fourier}
S_{l,m}=\sqrt{2/r} \, sin ({\pi}(l+1)(m+1)/r).
\end{equation}

While (\ref{twist}) and (\ref{fourier}) explicitly give the action
of a set of generators for $g=1$ ($\partial=\emptyset$), the action
in the higher genus case involves the $S$-matrix for the punctured
torus and $6j$-symbols. More specifically, $6j$-symbols (see
\cite[Chapter 7]{KL}) enter as the coefficients for the change of
basis corresponding to different labelled spines, see figure
\ref{6j} below.

Note that
${\rho}_k$ are {\em projective} representations of $SL(2,{\mathbb Z})$.
It is well known \cite{At} that they may be renormalized
to linear representations of $SL(2,{\mathbb Z})$ (and in fact of
$PSL(2,{\mathbb Z})$.)
More specifically, consider the presentation
$$SL(2,{\mathbb Z})\cong \; <\! S,T\, |\, S^4=1, S^2=(ST)^3 \! >.$$
The relations holds, under the image of
${\rho}_k$, only up to a phase. It is clear from the presentation that
the phases of ${\rho}_k(S)$, ${\rho}_k(T)$ may be adjusted so that
both relations actually hold in the image. A detailed analysis of the phase
factors \cite{MR} shows that they may be chosen to have a limit
as $k\longrightarrow \infty$.

\smallskip

\section{Kazhdan's property T} \label{property T}

We begin with the definition and a brief summary of implications of
property (T) (cf \cite{L} for a more detailed
discussion); the section concludes with the proof of theorem \ref{property T theorem}.
A discrete group $G$ has Kazhdan's property (T) if given
any (necessarily infinite dimensional)
unitary representation $\rho$ of $G$, whenever $\rho$ has almost
invariant vectors it also has a non-zero invariant vector. Recall
that a representation $\rho$ of a group $G$ {\em has almost
invariant vectors} if for any finite subset $S\subset G$ and for any
${\epsilon}>0$ there is a unit vector $v$ such that
$|{\rho}(s)v-v|<{\epsilon}$ for each $s\in S$.

Note that it suffices to check this condition for a fixed finite set
of generators of $G$: if it is true for $a$, $b$, then
$|{\rho}(ab)v-v|\leq |{\rho}(a)({\rho}(b)v-v)|
+|{\rho}(a)v-v|=|{\rho}(b)v-v|+|{\rho}(a)v-v|<2{\epsilon}$.

Property (T) may be expressed in terms of {\em Fell topology} on the set
$\widetilde G$ of unitary representations of $G$ (and similarly on the
set $\widehat G$ of irreducible unitary representations.) Let $H$ be a
Hilbert space, and ${\rho}\co G\longrightarrow {\cal U}(H)$ a
continuous unitary representation of $G$. Fell topology on $\widetilde
G$ is determined by open neighborhoods: let $S$ be a finite subset
of $G$, ${\epsilon}>0$ and $v\in H$ be a unit vector. Then the set
of representations ${\sigma}\co G\longrightarrow {\cal U}(H')$ for
which there exists a unit vector $v'\in H'$ such that
$ |(v,{\rho}(g)v)-(v',{\sigma}(g)v')|<{\epsilon} $
for each $g$ in $S$ is an open neighborhood of $(H,{\rho})$ in
$\widetilde G$. In these terms, $G$ has property (T) if and only if the
trivial one-dimensional representation ${\rho}_0$ is isolated in
$\widehat G$ with Fell topology, or equivalently if no
$\rho\neq {\rho}_0\in \widehat G$ {\em weakly contains} ${\rho}_0$.

Extending the meaning of the term {\em almost invariant} in the first
paragraph we call a sequence of vectors $\{ v_k \}$ belonging to
$\{ {\rho}_k \}$ {\em almost invariant} if for a generating set $S$,
$|{\rho}_k(s) v_k - v_k|\longrightarrow 0$ as $k\longrightarrow\infty$.
In our examples each ${\rho}_k$ is finite dimensional.

A short argument proves that infinite amenable groups don't have
property (T): suppose $G$ is amenable, and let $\{F_k \}$ be a
sequence of F\"{o}lner sets. Considering the characteristic
functions of these sets, it immediately follows that the left
regular representation of $G$ on $l^2(G)$ has almost invariant
vectors (and since the group is infinite it doesn't have a
non-trivial fixed vector.)

It follows that infinite abelian groups don't have property (T).
If $G$ has property (T) then so does any quotient group $Q$.
A consequence of these two statements is that
any group whose abelianization is infinite is not
(T). Another basic property is that a subgroup $H$ of finite index
in $G$ has (T) if and only if $G$ does. In particular,
$SL(2,{\mathbb Z})$ and the braid groups do not have property (T)
since they contain free non-abelian subgroups of finite index.

Finding a subgroup $H$ of finite index in a group $G$ with
$H^1(H;{\mathbb Z})$ infinite is a common way of proving that $G$ is
not (T). The question of existence of such subgroups in mapping
class groups is open in general. It is known that no such subgroups
$H$ exist in genus $g\geq 3$ with the additional assumption that $H$
contains the Torelli group \cite{Mc}, on the other hand such
subgroups exist in ${\cal M}_2$ \cite{T}, \cite{Mc}. This proves, in
particular, that ${\cal M}_2$ doesn't have property (T).

\noindent
{\em Proof of theorem} \ref{property T theorem}.
As above, given a surface $S$ of genus $g$ with $b$ boundary components and $n$
marked points, we denote its mapping class group by ${\cal M}^b_{g,n}$.
We have a central extension

$$1\longrightarrow {\mathbb Z}^n\longrightarrow {\cal M}^n_0\longrightarrow {\cal M}_{0,n}\longrightarrow 1$$

where the kernel is generated by Dehn twists along boundary parallel
curves.

Consider the hyperelliptic involution $J$ on the closed surface ${\Sigma}_g$
of genus $g$, inducing a double cover ${\Sigma}_g\longrightarrow S^2$
branched over $2g+2$ points. Let ${\cal M}^J_g$ denote the hyperelliptic
mapping class group of ${\Sigma}_g$: the subgroup consisting of elements
of ${\cal M}_g$ which commute
with $J$. Then we have a short exact sequence \cite{BH}

$$ 1\longrightarrow {\mathbb Z}/2{\mathbb Z}\longrightarrow
{\cal M}^J_g\longrightarrow {\cal M}_{0,2g+2}\longrightarrow 1$$

\noindent
where the generator of ${\mathbb Z}/2{\mathbb Z}$ is mapped to the
involution $J$. In particular, it suffices to show that
${\cal M}_{0,2g+2}$ doesn't have property (T): then the result for
${\cal M}^J_g$ is immediate.
Note that ${\cal M}^{J}_1={\cal M}_1\cong SL(2,{\mathbb Z})$
and ${\cal M}^{J}_2={\cal M}_2$ are known not to have property (T),
therefore our argument below specialized to $g=1,2$ gives another
proof for these groups.

Consider the family ${\rho}_k$ of projective representations  of
${\cal M}^n_0$, introduced in section \ref{quantum}. Fixing a spine,
a basis for the representations is indexed by labellings of its edges.
The strategy is to find almost invariant vectors $\{ v_k\}$ for $\{ {\rho}_k\}$
with {\em fixed boundary labels} and then to consider
${\rho}_k\otimes\bar{\rho}_k$
(the bar denotes the complex-conjugate representation.)
Then $v_k\otimes\bar{v}_k$ are almost invariant under
${\rho}_k\otimes\bar{\rho}_k$, and
since the extension ${\cal M}^n_0\longrightarrow {\cal M}_{0,n}$
is central, ${\rho}_k\otimes\bar{\rho}_k$ descends to a linear
representation of ${\cal M}_{0,n}$ (this follows from formula (\ref{twist})
for the action of the Dehn twists along the boundary parallel curves.)
Note that this also eliminates the projective ambiguity.

\begin{figure}[ht]
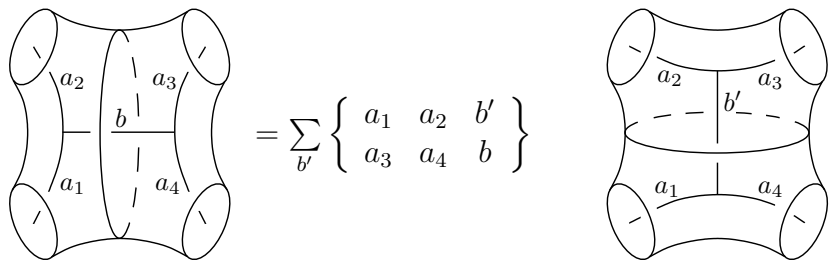

\sphere
\caption{Change of basis: the coefficients are given by $6j$ symbols.
The summation is over admissible
labels $b'$: in the context of the proof all labels $a_i$
equal $a=f(r )$ and $0\leq b'\leq 2a$.}
\label{6j}
\end{figure}

Consider the case of the $4$-punctured sphere ($n=4$), the general case
follows easily from this. The spine is H-shaped: there are four
``boundary'' labels $a_i$ and one ``interior'' label $b$, see figure \ref{6j}.
(Note that this is an {\em orthogonal} basis.)
Assume $r=k+2$ is large and fixed. Take all labels $a_i$ equal to
$a=f(r)$ where $f(r)$ is any function of $r$ with $f(r )\longrightarrow
\infty$ and $f(r)/{\sqrt r}
\longrightarrow 0$ as $r\longrightarrow\infty$.
Then it follows from the admissibility of the labels that $b$ ranges
from $0$ to $2a$ -- denote the corresponding basis vector by
$e_{a,b}$. (Here we omit reference to the level $k$ to simplify
the notation.)

Fix $a$ as above and consider the vector

$$v_k=\sum_{b=0}^{2a} \frac{e_{a,b}}{||e_{a,b}||}.$$

${\cal M}_{0,4}$ is generated by the Dehn twists ${\tau}_1,{\tau}_2$ along the
curves pictured in figure \ref{6j}:
one bounds a disk intersecting the bar of the H-spine in a point, the
other one corresponds to the dual I-spine. It follows from the choice of the
function $f(r)$ and formula (\ref{twist}) in section \ref{quantum}
that $v$ is almost invariant under
${\rho}_k({\tau}_1)$. The action of ${\rho}_k({\tau}_2)$ is easily understood
in the $I$-basis. Rewriting $v$ in this basis is hard since this requires
calculating the 6j-symbols (coefficients for the change of basis), however
due to the choice of the boundary labels and the admissibility (quantum
cut off) of labels, we know that the new interior label still satisfies $0\leq
b'\leq 2a$, making $v$ also almost invariant under ${\rho}_k({\tau}_2)$.

It remains to check that the unit vectors
$w_k:=v_k\otimes \bar{v}_k/||v_k\otimes \bar{v}_k||$ do not converge to
the fixed subspaces of ${\rho}_k\otimes\bar{\rho}_k$ as
$k\longrightarrow\infty$. According to formula (\ref{twist}), the basis vectors
$e_{a,b}\otimes \bar{e}_{a,c}$ are eigenvectors for ${\rho}_k({\tau}_1)$
with the corresponding eigenvalue

$$e^{2{\pi}i [b(b+2)-c(c+2)]/r}=e^{2{\pi}i (b-c)(b+c+2)/r}.$$

Therefore for prime $r$ the fixed vectors of
${\rho}_k\otimes\bar{\rho}_k$ have to lie in the subspace
spanned by $\{ e_{a,b}\otimes \bar{e}_{a,b},\, 0\leq b\leq r-2\}$ and by
$\{ e_{a,b}\otimes \bar{e}_{a,r-2-b}, \, 0\leq b\leq r-2\}$.

The component of the unit vector $w_k$ along the ``diagonal'': the span of
$\{ e_{a,b}\otimes \bar{e}_{a,b}\}$ tends to zero as $k\longrightarrow \infty$,
and since $b<2a\ll r$ the component along
$\{ e_{a,b}\otimes \bar{e}_{a,r-2-b}\}$ is trivial.
This concludes the proof for ${\cal M}_{0,4}$.

For ${\cal M}_{0,n}$, $n>4$, again choose a spine and fix all boundary labels
to be $a=f(r )$ as above. The tree structure of the spine (and in fact of
{\em all} spines) yields the required bound $b_j\ll \sqrt{r}$ for all
interior labels $b_j$ provided $r$ is large enough.
\qed

\smallskip

A similar argument allows one to find almost invariant vectors, with a fixed
boundary label, for representations ${\rho}_k$ of the mapping class group
${\cal M}_1^1$. Again ${\rho}_k\otimes\bar{\rho}_k$ descends to a representation
of ${\cal M}_1\cong SL(2,{\mathbb Z})$ with almost invariant vectors.
The situation is different for the {\em closed} torus: see section
\ref{tau} and the discussion in section \ref{comments}.

\bigskip

\section{The metaplectic (Segal-Shale-Weil) representation and the modular group} \label{metaplectic}

We first give a brief overview of the metaplectic representation,
defined on the double cover of the symplectic group $Sp(2n,{\mathbb
R})$ (cf \cite{Fo} for a detailed exposition.) The main discussion
will concern the case $n=1$, where restricting to the integers we
have a projective representation of $Sp(2,{\mathbb Z})\cong SL(2,
{\mathbb Z})$ which as mentioned in section \ref{quantum} can be
adjusted to be a linear representation of $SL(2,{ \mathbb Z})$. We
give an abstract definition in terms of intertwining operators and
then present an explicit formula in the $n=1$ case. The section
concludes with the proof of theorem \ref{convergence theorem}.

The Heisenberg group ${\mathbf H}_n$ is defined as ${\mathbb R}^{2n+1}$
with the coordinates $(p,q,t)=(p_1,\ldots,$ $p_n,q_1,\ldots,q_n,t)$ and with the
multiplication rule

$$ (p,q,t)(p',q',t')=(p+p',q+q',t+t'+\frac{1}{2}(pq'-qp')). $$

(This group is isomorphic to the group of $(n+2)\times (n+2)$ matrices
$M(p,q,t)$ with $1$'s on the diagonal
and with $M_{1,i}=p_{i-1}$, $M_{i,n+2}=q_{i+1}$, $i=2,\ldots,n+1$,
and $M_{n+2,n+2}=t$.)

Exponentiating the classical quantum-mechanical representation of
the Heisenberg algebra -- $P_j$ acts by $(1/2{\pi}i)\,\,
{\partial}/{\partial x_j}$, $Q_j$ acts by multiplication by $x_j$ --
one gets the Schr\"{o}dinger unitary representation ${\rho}$ of
${\mathbf H_n}$ on $L^2({\mathbf R}^n)$:

$${\rho}(p,q,t)f(x)=e^{2{\pi}it+2{\pi}iqx+{\pi}ipq} f(x+p). $$

The symplectic group $Sp(2n,{\mathbb R})$ acts on ${\mathbf H}_n$ by automorphisms:

$$ (p,q,t)\longmapsto (S(p,q),t) \;\, {\rm for} \;\, S\in Sp(2n,{\mathbb R}), $$

in particular
for any such $S$ we have a representation ${\rho}\circ S$ of ${\mathbf H_n}$ on
$L^2({\mathbb R}^n)$.

By the Stone-von Neumann theorem these representations are
equivalent: there exists a unitary operator ${\mu}(S)$ on
$L^2({\mathbb R}^n)$ such that ${\rho}\circ S={\mu}(S)\; {\rho}\;
{\mu}(S)^{-1}$. By Schur's lemma these intertwining operators define
a projective representation of the symplectic group:
${\mu}(TS)=C_{T,S}{\mu}(T){\mu}(S)$ for a scalar $C_{T,S}$ with
$|C_{T,S}|=1$; moreover it can be shown that it lifts to a linear
representation ${\mu}$ of the non-trivial double cover
$Mp(2n,{\mathbb R})$ of $Sp(2n,{\mathbb R})$.

The spaces $L^2_{\rm even}({\mathbb R}^n)$ and
$L^2_{\rm odd}({\mathbb R}^n)$ of even and odd functions
are invariant subspaces of the metaplectic representation $\mu$,
and it can be proved that the subrepresentations of $\mu$
on $L^2_{\rm even}$ and $L^2_{\rm odd}$ are irreducible
(and inequivalent).

In the case $n=1$ this representation has a particularly nice
explicit integral expression. Fix the standard generators of
$SL(2,{\mathbb Z})$:

$$ S=\left( \begin{array}{cc} 0 & 1\\ -1 & 0 \end{array} \right), \;
\;\; T=\left( \begin{array}{cc} 1&1\\0 & 1 \end{array} \right). $$

Then (up to a phase)
${\mu}(S)$ is the inverse of the Fourier transform and ${\mu}(T)$ is the
multiplication by the imaginary Gaussian:

$$ {\mu}(S)f(x)=\int_{\mathbb R} e^{2{\pi}ixy} f(y) \, dy, \; \; \;
{\mu}(T) f(x)=e^{-{\pi}ix^{2}} f(x). $$

It is easily seen that the subrepresentations
${\mu}'$, ${\mu}''$ of ${\mu}$ on
$L^2_{\rm odd}({\mathbb R})$ and $L^2_{\rm even}({\mathbb R})$
are given by
$${\mu}'(S)f(x)=\int_{\mathbb R} sin(2{\pi}xy) f(y) \, dy\;\;\;
{\rm and} \;\;\; {\mu}''(S)f(x)=\int_{\mathbb R} cos(2{\pi}xy) f(y) \, dy $$

(in both cases $T$ acts by multiplication by $e^{-{\pi}ix^{2}}$.)

Theorem \ref{convergence theorem} asserts that the representations
${\rho}_k$ converge in Fell topology to ${\mu}'$. For the proof,
consider the formulas (\ref{twist}) and (\ref{fourier}) in section
\ref{quantum} for the action of $T$ and $S$.
It is convenient to rescale the variables: $a=l/{\sqrt{2r}}$,
$b=m/{\sqrt{2r}}$. We always use the notation $r=k+2$
and in our estimates we assume $r$ is large.
The new parameters range from $0$ to $\sim{\sqrt{r/2}}$ and values are spaced
$1/{\sqrt{2r}}$ apart.
Denoting the basis vector for ${\rho}_k$ corresponding to the
label $l$ (and to the rescaled value $a$) by $e_a$,
the representation is given by

$$ {\rho}_k(T)e_a\simeq e^{-{\pi}i a^2}e_a, \;\;\;
{\rho}_k(S)e_a\simeq \sum_b\sqrt{2/r} \, sin \, (2{\pi}ab) e_b
$$

According to the definition of Fell topology, given an odd function
$f\in L^2({\mathbb R})$ of norm $1$ we need to find a sequence of
unit vectors $\{ v_k\}$, $v_k\in {\rho}_k$ such that $\forall {\epsilon}>0$

$$|(f,{\mu}'(S)f)_{L_2}-(v_k,{\rho}_k(S) v_k)|<{\epsilon} \;\;\; {\rm and}
\;\;\; |(f,{\mu}'(T)f)_{L_2}-(v_k,{\rho}_k(T) v_k)|<{\epsilon}$$

for sufficiently large $k$.
We will view vectors in ${\rho}_k$ as step functions on
$[0,{\sqrt{r/2}}]$, with the step length equal to $1/\sqrt{2r}$.
Since step functions are dense in
$L^2({\mathbb R})$, for a sufficiently large $k$ there exists
a vector $v_k\in {\rho}_k$ whose corresponding step function $g$ satisfies
$|g-f|_{L^2}<{\epsilon}/3$. More precisely, note that the given function
$f$ is {\em odd} while the available step functions are supported on the
interval $[0, {\sqrt{r/2}}]$ which is in ${\mathbb R}_{+}$.
Therefore the correspondence between step functions $g$ and vectors
$v_k$ involves keeping only the $x\geq 0$ half of $g$ and multiplying
by $\sqrt 2$ to keep the norm equal to $1$.
Increasing $k$ if necessary we may assume also that

$$|(g,{\mu}'(S)g)_{L^2}-(v_k,{\rho}_k(S) v_k)|<{\epsilon}/3.$$
Then
$$|(f,{\mu}'(S)f)-(v_k,{\rho}_k(S) v_k)|\leq
|(f,{\mu}'(S)f)-(g,{\mu}'(S)g)|+{\epsilon}/3 \leq$$
$$|(f-g,{\mu}'(S)f)|+
|(g,{\mu}'(S)(f-g)|+{\epsilon}/3\leq {\epsilon}.$$

Here we used the fact that ${\mu}'(S)$ is an isometry.
The estimate for the action of $T$ is proved identically.
\qed

\bigskip

\section{Property $\tau$, the metaplectic representation over
${\mathbb Z}/p{\mathbb Z}$, and an analytic fact.} \label{tau}

The question of whether the quantum representations ${\rho}_k$ of
the mapping class group of the closed torus -- $SL(2,{\mathbb Z})$ --
have almost invariant vectors leads to an interesting
problem in Fourier analysis. In particular, the sequence of
Gaussians $e^{-{\pi}l^2/r}$ ``almost'' solves the problem
as $r\longrightarrow\infty$: they are almost invariant under
the Fourier transform, however their support is exactly
balanced against the quadratic behavior
of the eigenvalues (\ref{twist}) of the Dehn twist $T$.

In the limit the problem is to find a function in $L^2({\mathbb R})$
almost invariant under the Fourier transform $\;\,\widehat{   }\; \,$
(more precisely, under
the $sin$-transform) and supported near the set of points
$\{ \sqrt{n},\, n\in{\mathbb Z} \}$.
It follows from the uncertainty principle that a solution cannot
be given by a function similar to the Gaussian:
there doesn't exist a function supported near zero,
say on $(-{\epsilon},{\epsilon})$ -- making it almost invariant
under the action of $T$ -- whose Fourier transform is also mostly
supported on $(-{\epsilon},{\epsilon})$ for small $\epsilon$.

In fact, as stated in lemma \ref{factor} in the introduction,
the quantum representations
${\rho}_k$ (for $k$ prime) factor through $SL(2,{\mathbb Z}/k{\mathbb Z})$
(the proof of this is given below.) It follows from Selberg's theorem
\cite{L} that $SL(2,{\mathbb Z})$ has property
$\tau$ with respect to the family of congruence subgroups
$$ N_k=ker[SL(2,{\mathbb Z})\longrightarrow SL(2,{\mathbb Z}/k{\mathbb Z})],$$

in other words the trivial representation is isolated in the family
of representations $\{ {\rho}\, |\, ker({\rho})\supset N_k$ for some $k\}$.
This proves corollary \ref{no almost invariant vectors}:
our representations ${\rho}_k$ of $SL(2,{\mathbb Z})$ do not have
almost invariant vectors; in particular the problem concerning
the Fourier transform on $\mathbb R$ above does not have a solution.
In particular we can prove:

\medskip

\begin{thm} \label{analytic}
\sl There are sets $X\subset {\mathbb R}$ of constant density
so that no set similar to $X$ supports almost invariant functions.
\end{thm}

\medskip

\noindent
{\bf Notations}: We say $X$ has {\em constant density} if there exist
$N,c>0$ such that

$$ {\mu}(X\cap [r,r+N])/N>c\,\,\, {\rm for}\,\, {\rm all}\,\,\,
r\in {\mathbb R} $$

where $\mu$ is the Lebesgue measure. We say $X$ supports {\em almost invariant
functions} if there exist $L^2$ functions $f_i\co X\longrightarrow
{\mathbb C}$ such that $||f_i- \widehat{f_i}||\longrightarrow 0$.
\hspace{.05cm} $\widehat{   }$ \hspace{.1cm} denotes
the Fourier transform.
We say $Y$ is {\em similar} to $X$ if there exist $r_1\in
{\mathbb R}\smallsetminus 0$ and $r_2\in {\mathbb R}$ such that
$Y=r_1X+r_2$.

\noindent
{\bf Proof.} Let, for ${\epsilon}>0$, $X_{\epsilon}=\{ r\in {\mathbb R}|
\, |r-(\pm\sqrt{n})|<{\epsilon}/\sqrt{n}, n$ a positive integer$\}$.
Clearly all $X_{\epsilon}$ have constant density. If each supported
almost invariant functions, a diagonalization argument and the
preceding discussion would contradict property $\tau$ for
$SL(2,{\mathbb Z})$. If $X$ is replaced by $r_1 X$ rescale
$a=l/r_1 \sqrt{r}$, $b=m/r_1\sqrt{r}$; translation by $r_2$ introduces a phase
$e^{2{\pi}ir_2}$ into $\hat{f}$ but since our analysis is on ${\rho}_k\otimes
\bar{{\rho}_k}$ such phases are irrelevant.
\qed.

\smallskip

To prove lemma \ref{factor} we note that the construction of the
metaplectic representation $\mu$, sketched in the previous section,
can be carried through for finite fields $K$ in place of ${\mathbb R}$,
see \cite{N} and references therein. As in the case of $\mathbb R$,
one considers the action of the symplectic group $Sp(2n,K)$ on the
Heisenberg group ${\mathbf H_n}$ (defined in terms of $K$), and the
Stone-von Neumann theorem asserts that the Schr\"{o}dinger representation
is the unique unitary representation ${\mathbf H_n}$ with the given
action of the center of ${\mathbf H_n}$.
Again, the irreducible subspaces consist of even and odd functions.

Theorem 4.1 in \cite{N} gives an explicit formula for the metaplectic
representation
$\mu$ over finite fields, in particular taking $K={\mathbb Z}/p{\mathbb Z}$
one has (for $n=1$ and the standard generators $S$, $T$ of
$Sp(2,{\mathbb Z}/p{\mathbb Z})\cong SL(2,{\mathbb Z}/p{\mathbb Z})$)

\begin{equation} \label{finite}
{\mu}(T) f(x)= e^{-{\pi}i x^2/2p} f(x), \; \; \; \;
{\mu} (S) f(x) = \frac{1}{\sqrt{p}} \sum_{y\in K}
e^{{\pi}ixy/p} f(y)
\end{equation}

\noindent
Comparing (\ref{finite}) with the formulas (\ref{twist}), (\ref{fourier}) in
section \ref{quantum} and restricting to the subspace of odd functions,
observe that (for prime $p$) the quantum representation ${\rho}_p$
may be viewed as the irreducible component ${\mu}'$ of the metaplectic
representation of $SL(2,{\mathbb Z}/p{\mathbb Z})$.
\qed

\bigskip

\section{Comments and questions} \label{comments}

{\bf 1. Property (T) and mapping class groups}.
We established in section \ref{tau} that the family of quantum
$SU(2)$-representations of $SL(2,{\mathbb Z})$ that arise as the
mapping class group of the closed torus doesn't have almost invariant
vectors. On the other hand, the representations of the mapping class groups
of the punctured torus, and of the sphere with $4$ (or more) punctures do have
almost invariant vectors. It is an interesting question which of the two
possibilities
holds for closed surfaces of higher genus. (Of course the latter one would
imply that the mapping class groups of all orientable surfaces do not have
property (T).)
An answer to this question would likely require a detailed analysis
of the action of a generating set of Dehn twists on explicit vectors under
the quantum representations.
Checking that a proposed vector $v$ is almost invariant would involve
calculations with $6j$-symbols and with the $S$-matrix for
the punctured torus. The last ingredient is known \cite{Ki} in terms of
Macdonald polynomials but the calculations are difficult in practice.

{\bf 2. Limits of quantum representations.}
We showed that the representations for the closed torus converge,
as the level approaches infinity, to the metaplectic representation
of $SL(2, {\mathbb Z})$.
It would be interesting to find out whether there is a limiting
representation of the mapping class group ${\cal M}_g$ in the higher genus case.

Specifically consider the higher-dimensional metaplectic representation
${\rho}$. The maps induced on the first homology of the surface by
homeomorphisms give a surjection
${\pi}\co {\cal M}_g\longrightarrow Sp(2g,{\mathbb Z})$. The quantum $SU(2)$ representations
faithfully detect the mapping class groups modulo center \cite{A}, \cite{FWW}
so they cannot factor through the symplectic group.
However it is possible that the quantum representations ${\rho}_k$ converge
to ${\rho}\circ{\pi}$ in Fell topology. Convergence in this topology is rather
weak and proving this would require exhibiting, for any function
in $L^2({\mathbb R}^{2g})$, a sequence of vectors $v_k$ that are transformed
by ${\rho}_k$ applied to a generating set of Dehn twists approximately as
by ${\rho}\circ{\pi}$ -- generalizing our proof in the genus $1$ case.

\end{document}